\documentclass{article}

\usepackage{amssymb}
\usepackage{amsmath}
\usepackage[all]{xy}
\usepackage{a4wide}
\usepackage{enumerate,fancyhdr}

\begin{document}

\title{Classe d'isog\'enie de vari\'et\'es ab\'eliennes pleinement de type GSp}

\date{}

\author{Nicolas Ratazzi \footnote{Universit\'e Paris-Sud XI, B\^atiment 425, 91405 Orsay Cedex, FRANCE, nicolas.ratazzi@math.u-psud.fr}}

\renewcommand{\epsilon}{\varepsilon}
\newcommand{\Gal}{\textnormal{Gal}}
\newcommand{\Z}{\mathbb{Z}}
\newcommand{\Q}{\mathbb{Q}}
\newcommand{\C}{\mathbb{C}}
\newcommand{\R}{\mathbb{R}}
\newcommand{\N}{\mathbb{N}}
\newcommand{\F}{\mathbb{F}}
\newcommand{\G}{\mathbb{G}}
\newcommand{\T}{\textnormal{T}}
\newcommand{\Frob}{\textnormal{Frob}}
\newcommand{\Id}{\textnormal{Id}}
\newcommand{\PSp}{\textnormal{PSp}}
\renewcommand{\O}{\mathcal{O}}
\newcommand{\End}{\textnormal{End}}
\renewcommand{\H}{\textnormal{H}}
\newcommand{\MT}{\textnormal{MT}}
\newcommand{\GL}{\textnormal{GL}}
\newcommand{\SL}{\textnormal{SL}}
\newcommand{\Sp}{\textnormal{Sp}}
\newcommand{\PGSp}{\textnormal{PGSp}}
\newcommand{\GSp}{\textnormal{GSp}}
\renewcommand{\text}{\textnormal}

\newtheorem{theo}{Th{\'e}or{\`e}me} [section]
\newtheorem{lemme}[theo]{Lemme}
\newtheorem{conj}{Conjecture}[section]
\newtheorem{prop}[theo]{Proposition}
\newtheorem{cor}[theo]{Corollaire}
\newtheorem{defi}[theo]{D\'efinition}
\newtheorem{notations}{Notations}
\newtheorem{rem}[theo]{Remarque}

\newcommand{\demo}{\noindent \textit{D{\'e}monstration} : }
\newcommand{\findemo}{\hfill $\Box$}

\maketitle

\noindent \textbf{Abstract : } Faltings in 1983 proved that a necessary and sufficient condition for two abelian varieties $A$ and $B$ to be isogenous over a number field $K$ is that the local factors of the L-series of $A$ and $B$ are equal for almost all primes of $K$ ; for each such prime this implies that $A$ and $B$ have the same number of points over the residue field. We show in this article that for abelian varieties faithfully of type $\GSp$ (a class containing the abelian varieties with endomorphism ring $\Z$ and of odd dimension) `having the same number of points' may be replaced by `the number of points have the same prime divisors' and still gives a sufficient condition for $A$ and $B$ to be $K$-isogenous. The proof is based on ideas of Serre \cite{serreim72} and Frey-Jarden \cite{FJ} and follows closely Hall-Perucca \cite{hallp} who proved the result for elliptic curves.

\section{Introduction}

\noindent En 1983, Faltings a prouv\'e le r\'esultat suivant pour deux vari\'et\'es ab\'eliennes $A$ et $B$ d\'efinies sur un corps de nombres~: une condition n\'ecessaire et suffisante pour que $A$ et $B$ soient $K$-isog\`enes est qu'il existe un ensemble $S$ de densit\'e un de premiers $\mathfrak{p}$ de $K$ (de bonne r\'eduction pour $A$ et $B$) tels que les facteurs locaux des s\'eries L de $A$ et $B$ soit \'egaux. Ce r\'esultat a \'et\'e r\'ecemment am\'elior\'e dans \cite{hallp} pour les courbes elliptiques de la fa\c{c}on suivante : plut\^ot que demander que les r\'eductions $A_{\mathfrak{p}}$ et $B_{\mathfrak{p}}$ aient le m\^eme nombre de points sur le corps r\'esiduel $k_{\mathfrak{p}}$  (condition \'equivalente \`a l'\'egalit\'e des facteurs locaux des s\'eries L dans le cas de dimension 1), il suffit (pour des courbes elliptiques) de demander que le nombre de points de $A(k_{\mathfrak{p}})$ et de $B(k_{\mathfrak{p}})$ aient le m\^eme ensemble de diviseurs premiers. De plus il suffit de savoir ceci non pour tous les nombres premiers mais seulement pour une famille infinie. Dans cet article nous am\'eliorons le r\'esultat de Faltings de la m\^eme fa\c{c}on pour les vari\'et\'es ab\'eliennes pleinement de type GSp (au sens de la d\'efinition \ref{pleinement} ci-dessous), famille contenant les vari\'et\'es ab\'eliennes ayant un anneau d'endomorphismes $\Z$ et de dimension $2$ ou impaire. Nous suivons pour cela la strat\'egie de Hall-Perucca \cite{hallp}, elle m\^eme bas\'ee sur des travaux ant\'erieurs de Serre \cite{serreim72} et de Frey-Jarden \cite{FJ}.

\subsection{Vari\'et\'es ab\'eliennes pleinement de type GSp\label{va}}

\noindent Soit $A$ une vari\'et\'e ab\'elienne de dimension $g\geq 1$ d\'efinie sur un corps de nombres $K$. Consid\'erant l'action du groupe de Galois $G_K :=\Gal(\bar K/K)$ sur les points de $\ell^{\infty}$-torsion, pour $\ell$ premier, on associe naturellement \`a $A/K$, la repr\'esentation $\ell$-adique 
\[\rho_{\ell^{\infty},A}: G_K \rightarrow\GL(\T_{\ell}(A))\simeq\GL_{2g}(\Z_{\ell})\]

\noindent avec $\T_{\ell}(A)=\varprojlim A[\ell^n]$ le module de Tate $\ell$-adique de $A$ et on note \'egalement son image

\[G_{\ell^{\infty},A}:=\rho_{\ell^{\infty},A}\left(G_K\right).\]

\noindent Nous noterons \'egalement $\rho_{\ell,A}$ et $G_{\ell,A}$  (voire $\rho_{\ell}$  et $G_{\ell}$ s'il n'y a pas d'ambiguit\'e) les objets d\'eduits modulo $\ell$. On suppose fix\'ee une polarisation sur $A$. Dans ce cas les $T_{\ell}(A)$ sont munis, pour $\ell$ premier au degr\'e de la polarisation, d'une forme altern\'ee $e_{\ell}$ et la repr\'esentation  $\rho_{\ell^{\infty},A}$ est \`a valeurs dans $\GSp_{2g}(T_{\ell}(A),e_{\ell})\simeq  \GSp_{2g}(\Z_{\ell})$.

\medskip

\noindent Un probl\`eme naturel est de savoir quand l'image $G_{\ell^{\infty},A}$ est d'indice fini dans (voire \'egale \`a) $\GSp_{2g}(\Z_{\ell}),$ pour tout premier $\ell$ assez grand (d\'ependant de $A/K$ et premier au degr\'e de la polari\-sation).

\begin{defi} \label{pleinement} Soit $A/K$ une vari\'et\'e ab\'elienne polaris\'ee de dimension $g\geq 1$. Nous dirons que $A$ est \textit{pleinement de type $\GSp$} si $A$ est telle que 
\[\textit{pour tout premier $\ell$ assez grand},\ \ G_{\ell,A}=\GSp_{2g}(\F_{\ell}).\]
\end{defi}

\noindent On sait par exemple apr\`es Serre \cite{serreim72} que toute courbe elliptique sans CM, ie \`a anneau d'endo\-morphismes (sur $\bar K$) $\End_{\bar K}(A)=\Z$, est pleinement de type $\GSp$ (il s'agit m\^eme dans ce cas d'une condition \'equivalente \`a \^etre sans CM sur $\bar K$). En dimension quelconque, une condition n\'ecessaire est d'avoir $\End_{\bar K}A=\Z$~; cette condition n'est pas en g\'en\'eral suffisante, mais on sait qu'elle l'est (cf. th\'eor\`eme \ref{th1} ci-dessous) si $g$ n'appartient pas \`a l'ensemble exceptionnel $\Sigma$ d\'efini comme suit~:
\begin{equation}\label{defdeS}
\Sigma=\left\{g\geq 1\;|\; \exists k\geq 3,\;{\rm impair},\; \exists a\geq 1,\; 2g=(2a)^{k} \text{ ou } 2g={2k\choose k}\right\}.
\end{equation}

\noindent Dans notre contexte le r\'esultat important est  un  th\'eor\`eme de Serre (\cite{college8586} Th\'eor\`eme 3 et \cite{serremfv} Th\'eor\`eme 3) compl\'et\'e par Pink (\cite{pink1} Theorem 5.14) et dans une autre direction par Hall (\cite{hall} Theorem 1).

\begin{theo}\label{th1}\textnormal{\textbf{(Serre, Pink, Hall)}} Si $A/K$ est une vari\'et\'e ab\'elienne de dimension $g$ n'appar\-tenant pas \`a $\Sigma$, d\'efinie sur un corps de nombres, telle que $\End_{\bar K}(A)=\Z$, alors $A$ est pleinement de type $\GSp$. Si $g$ est quelconque mais l'on suppose que le groupe de Mumford-Tate $\MT(A)=\GSp$ et que le mod\`ele de N\'eron de $A$ sur l'anneau des entiers $\mathcal{O}_K$ poss\`ede une fibre semistable avec dimension torique \'egale \`a un,  la m\^eme conclusion vaut.
\end{theo}

\subsection{R\'esultat principal\label{structure}}

\noindent Le th\'eor\`eme \ref{th1} pr\'ec\'edent donne une vaste classe de vari\'et\'es ab\'eliennes pleinement de type $\GSp$. Nous pouvons maintenant \'enoncer notre r\'esultat principal. 

\begin{theo}\label{radical} Soit $K$ un corps de nombres et soient $A_1,A_2/K$ deux vari\'et\'es ab\'eliennes pleinement de type $\GSp$. Consid\'erons $S$ un sous-ensemble de places finies $v$ de $K$, de corps r\'esiduel $\F_v$, de bonne r\'eduction pour $A_1$ et $A_2$, de densit\'e analytique $1$ et supposons \'egalement donn\'e un sous-ensemble $\Lambda$ infini de l'ensemble des nombres premiers. Alors $A_1$ est $K$-isog\`ene \`a $A_2$ si et seulement si
\[\forall v\in S,\ \forall \ell\in\Lambda\ \ \left( \ell\mid \text{Card}(A_1(\F_v))\iff \ell\mid \text{Card}(A_2(\F_v))\right).\]
\end{theo}

\noindent Notons que l'on sait apr\`es Faltings \cite{falt} que la classe d'isog\'enie d'une vari\'et\'e ab\'elienne $A/K$ est donn\'ee par sa fonction $\zeta$ qui donne en particulier les valeurs $\text{Card}(A(\F_v))$. Le th\'eor\`eme \ref{radical} (dont la preuve utilise n\'eanmoins \cite{falt}) prouve qu'une donn\'ee sensiblement plus faible est en fait suffisante (au moins pour les vari\'et\'es ab\'eliennes pleinement de type $\GSp$).

\medskip

\noindent \textbf{Structure de la preuve : } Nous commen\c{c}ons par des rappels sur les groupes symplectiques utilis\'es dans les parties suivantes. Les paragraphes 3 et 4 sont consacr\'es \`a la preuve proprement dite du th\'eor\`eme \ref{radical}. Pour cela nous suivons de pr\`es en l'adaptant en dimension sup\'erieure l'argument de Hall-Perucca \cite{hallp} et les id\'ees de Serre \cite{serreim72} et Frey-Jarden \cite{FJ}. La structure de la preuve de l'implication \textit{si} (l'autre implication \'etant facile) du th\'eor\`eme \ref{radical} est en trois \'etapes : 
\begin{enumerate}
\item \label{1} Montrer, pour $\ell$ variant dans un sous-ensemble infini $\Lambda_1$ de $\Lambda$, que $K(A_1[\ell])=K(A_2[\ell])$, puis en d\'eduire que $A_1$ et $A_2$ sont $\bar{K}$-isog\`enes.
\item \label{2} Montrer qu'il existe un caract\`ere quadratique $\epsilon : G_K\rightarrow \{\pm 1\}$ tel que pour tout $\ell$ variant dans un sous-ensemble infini $\Lambda_2$ de $\Lambda_1$ les repr\'esentations $\rho_{\ell,A_1}$ et $\epsilon\otimes\rho_{\ell,A_2}$ sont isomorphes.
\item \label{3} Montrer que pour tout $\ell$ variant dans un sous-ensemble infini $\Lambda_3$ de $\Lambda_2$, la repr\'esentation $\rho_{\ell,A_1}$ est en fait isomorphe \`a $\rho_{\ell,A_2}$, puis conclure que $A_1$ et $A_2$ sont $K$-isog\`enes.
\end{enumerate}

\medskip

\noindent Pour l'\'etape 1 nous avons besoin de prouver un r\'esultat d'un int\'er\^et ind\'ependant : le th\'eor\`eme d'isog\'enies horizontales \ref{isogenies1} (suivant la terminologie de \cite{FJ}) pour la classe de vari\'et\'es ab\'eliennes consid\'er\'ees. Ce th\'eor\`eme nous assure que deux vari\'et\'es ab\'eliennes pleinement de type GSp sont isog\`enes sur $\bar K$ d\`es qu'elles v\'erifient une condition du type $\left[K\left(A_1[\ell],A_2[\ell]\right):K\left(A_1[\ell]\right)\right]\leq c$ pour une certaine constante $c>0$ et pour une infinit\'e de premiers $\ell$. Ce r\'esultat peut se voir comme un corollaire du th\'eor\`eme 1.6 de \cite{hrJIMJ} prouvant la conjecture de Mumford-Tate forte pour le produit $A_1\times A_2$ de telles vari\'et\'es ab\'eliennes. Il remplace dans notre preuve l'usage fait par \cite{hallp} du Theorem A de \cite{FJ}. Par ailleurs nous aurons besoin d'un raffinement de ce r\'esultat et nous donnons donc une preuve alternative, dans l'esprit de \cite{FJ} et \cite{serreim72}, notamment en prouvant le lemme \ref{reduction}.

\medskip

\noindent La proposition \ref{twist} donne essentiellement une preuve de l'\'etape 2. Il s'agit d'une adaptation en dimension sup\'erieure du paragraphe 6.2 de \cite{serreim72}, en utilisant les r\'esultats de Raynaud \cite{raynaud} en lieu et place du paragraphe 1 de \cite{serreim72}.

\medskip

\noindent L'\'etape 3 se prouve comme dans \cite{hallp} et utilise notamment le th\'eor\`eme de Faltings (proposition \ref{lemfalt} rappel\'ee plus loin).

\medskip

\noindent \textbf{Remerciements~:} Je remercie Pascal Autissier, Chris Hall et Marc Hindry pour les commentaires qu'ils m'ont fait sur une version pr\'eliminaire de ce texte. Je remercie \'egalement le rapporteur dont les suggestions m'ont permis d'am\'eliorer la pr\'esentation de la partie 3 de cet article.

\section{Rappels sur les groupes symplectiques}

\noindent \'Etant donn\'e un entier $g\geq 1$, on rappelle la d\'efinition du \textit{groupe (alg\'ebrique) symplectique}~:

\[\GSp_{2g}:=\left\{M\in\GL_{2g}\;|\; \exists \lambda(M)\in\G_m,\;^t\!{M}\begin{pmatrix} 0 & I_g\cr -I_g&0\cr\end{pmatrix}M=\lambda(M) \begin{pmatrix} 0 & I_g\cr -I_g&0\cr\end{pmatrix}\right\}.\]

\noindent C'est un groupe alg\'ebrique  sur $\Z$. On introduit $\lambda~:~\GSp_{2g}\rightarrow\G_m$, l'homomorphisme qui associe \`a $M$ son multiplicateur $\lambda(M)$. Par d\'efinition, le \textit{groupe sp\'ecial symplectique} $\Sp_{2g}$ est le noyau du morphisme $\lambda$ et un corps $F$ \'etant donn\'e, les groupes \textit{projectifs symplectiques} $\PGSp_{2g}(F)$ et $\PSp_{2g}(F)$ sont les quotients respectifs de $\GSp_{2g}(F)$ et de $\Sp_{2g}(F)$ par les matrices scalaires appartenant respectivement \`a $\GSp_{2g}(F)$ et \`a $\Sp_{2g}(F)$.

\begin{rem}\label{detlambda} Notons le lien suivant entre l'application multiplicateur et le d\'eterminant~:
\[\forall M\in\GSp_{2g},\ \ \  \det M=\lambda(M)^g.\]
\end{rem}

\noindent Rappelons \'egalement le lien bien connu entre caract\`ere cyclotomique, repr\'esentation $\ell$-adique et multiplicateur : 

\begin{lemme}\label{cyclo} Soit $\ell$ un premier ne divisant pas le degr\'e de la polarisation fix\'ee sur une vari\'et\'e ab\'elienne polaris\'ee $A$. En notant $\rho_{\ell}$ la repr\'esentation $\ell$-adique modulo $\ell$ associ\'ee \`a $A/K$, la compos\'ee $\lambda\circ \rho_{\ell} : \Gal(\bar{K}/K)\rightarrow \F_{\ell}^{\times}$ n'est autre que le caract\`ere cyclotomique $\chi_{\text{cycl}}$. 
\end{lemme}

\begin{lemme}\label{cardinal} Soient $g_1, g_2$ deux entiers strictement positifs et soit $\ell$ un nombre premier impair. On se donne \'egalement un sous-groupe $Z_i\subset \{\pm I_{2g_{i}}\}\subset \Sp_{2g_i}(\F_{\ell})$ pour $i\in\{1,2\}$. On a alors
\[ \left|\Sp_{2g_1}(\F_{\ell})/Z_1\right|=\left|\Sp_{2g_2}(\F_{\ell})/Z_2\right|\Rightarrow g_1=g_2 \text{ et } Z_1=Z_2.\]
\end{lemme}
\demo Le cardinal d'un $\ell$-Sylow est $\ell^{g_i^2}$ donc $g_1=g_2$ et on en tire imm\'ediatement la seconde assertion. \findemo

\medskip

\noindent Nous rappelons ici un lemme tir\'e de \cite{hrJIMJ} (cf. lemmes 2.13 et 2.15 ainsi que les preuves), rassemblant quelques r\'esultats classiques sur les groupes symplectiques, leurs sous-groupes distingu\'es et leurs automorphismes (voir notamment Dieudonn\'e~\cite{dieu} Chap. IV paragraphes 3 et 6).

\begin{lemme}\label{gsimple} Soit $g\geq 1$ et soit $\F$ un corps~; on exclut les cas $g=1$, $\F=\F_2$, $\F_3$ ou $\F_4$ et $g=2$, $\F=\F_2$.
\begin{enumerate}
\item Le seul sous-groupe distingu\'e non trivial de $\Sp_{2g}(\F)$ est son centre $\{\pm I_{2g}\}$.
\item Les $\F$-automorphismes de $\PGSp_{2g}(\F)$ sont tous int\'erieurs et l'action par conjugaison se fait via un \'el\'ement de $\Sp_{2g}(\F)$.
\end{enumerate}
\end{lemme} 

\section{Isog\'enies horizontales\label{MTS}}

\noindent Nous voulons ici prouver le th\'eor\`eme \ref{isogenies1} ci-dessous d'isog\'enies horizontales (dans la terminologie de \cite{FJ}). Nous donnons d'abord une preuve directe utilisant le th\'eor\`eme 1.6 de \cite{hrJIMJ} prouvant la conjecture de Mumford-Tate forte pour un produit de vari\'et\'es ab\'eliennes pleinement de type $\GSp$. Puis nous donnons, dans la suite de ce paragraphe une preuve alternative, (not\'ee preuve (B) ci-apr\`es) utilisant directement les id\'ees du paragraphe 6 de Serre \cite{serreim72}. La raison est que les divers \'enonc\'es qui interviennent dans cette preuve (B), notamment la proposition \ref{twist}, sont \'egalement n\'ecessaires pour la preuve du th\'eor\`eme principal \ref{radical}.

\medskip 

\noindent Dans la suite $K$ est un corps de nombres et $A_1$ et $A_2$ des vari\'et\'es ab\'eliennes pleinement de type $\GSp$, de dimensions respectives $g_1$, $g_2$, d\'efinies sur $K$. La notation $\ell$ indiquera toujours un nombre premier. 

\medskip

\noindent \textbf{Notation} Dans la suite, nous utiliserons la notation $\gg$ pour signifier \textquotedblleft sup\'erieur \`a, \`a une constante multiplicative pr\`es (d\'ependant \'eventuellement de $A_1$, $A_2$ mais) ind\'ependante de $\ell$\textquotedblright . De m\^eme pour la notation $\ll$~; la notation $\gg\ll$ signifiant la conjonction de $\gg$ et de $\ll$. 

\medskip

\noindent Pour $i\in\{1,2\}$, posons $K_{\ell,i}:=K(A_i[\ell])$. Par d\'efinition de pleinement de type $\GSp$, on a, pour $\ell\gg 1$, les diagrammes d'extensions suivants :
\[
	\xymatrix{
		K_{\ell,i}\ar@{-}[dd]_{\GSp_{2g_i}(\F_\ell)} & \\
		   & K(\mu_\ell)  \ar@{-}[ul]_{\Sp_{2g_i}(\F_\ell)} \\
		K  \ar@{-}[ur]_{\F_\ell^\times}
	}
\]

\noindent Notons par ailleurs $K_{\ell}^{\cup}$ la $K$-extension engendr\'ee par $K_{\ell,1}\cup K_{\ell,2}$  et $\ K_{\ell}^{\cap}$ la $K$-extension (engendr\'ee par) $K_{\ell,1}\cap K_{\ell,2}$. Enfin, pour $i\in\{1,2\}$, posons $H_{\ell,i}:=\Gal(K_{\ell,i}/K_{\ell}^{\cap})$. Par d\'efinition de pleinement de type $\GSp$, on a, pour $\ell\gg 1$, le diagramme d'extensions suivant :

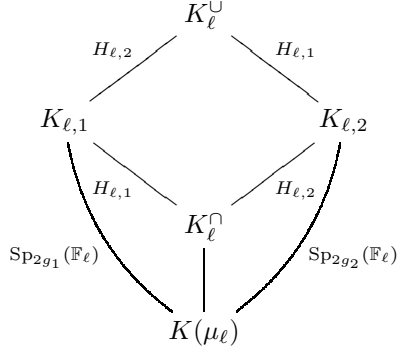
\[
	\xymatrix{
		& K_{\ell}^{\cup} & \\
	K_{\ell,1} \ar@{-}[ur]^{H_{\ell,2}} & & K_{\ell,2} \ar@{-}[ul]_{H_{\ell,1}} \\
	    & K_{\ell}^{\cap} \ar@{-}[ur]_{H_{\ell,2}} \ar@{-}[ul]^{H_{\ell,1}} & \\
	    & K(\mu_\ell)
	    	\ar@{-}[u]
			\ar@{-}@/_{1pc}/[uur]_{\Sp_{2g_2}(\F_\ell)}
	    	\ar@{-}@/^{1pc}/[uul]^{\Sp_{2g_1}(\F_\ell)} &
	}
\]

\begin{lemme}\label{reduction0} Supposons qu'il existe $c>0$ et un ensemble $\Lambda$ infini tel que pour tout $\ell\in \Lambda$, on a 
\[ \text{Card}(H_{\ell,2})\leq c.\]
\noindent Alors, pour $i\in\{1,2\}$ et $\ell\in\Lambda$, on a 
\[\ell\gg 1\Rightarrow \text{Card}(H_{\ell,i})\leq 2.\]
\end{lemme}
\demo Consid\'erons le diagramme d'extensions pr\'ec\'edent. L'extension $K_{\ell}^{\cap}/K(\mu_{\ell})$ est galoisienne, donc si $\ell$ est assez grand, les vari\'et\'es ab\'eliennes $A_i$ \'etant pleinement de type $\GSp$, les groupes $H_{\ell,i}$ sont des sous-groupes distingu\'es de $\Sp_{2g_i}(\F_{\ell})$ pour $i\in\{1,2\}$. Par le lemme \ref{gsimple} on en d\'eduit que ces groupes sont soit $\{I_{2g_i}\}$, soit $\{\pm I_{2g_i}\}$ (pour $i\in\{1,2\}$), ou bien le groupe sp\'ecial symplectique tout entier.  Or si $H_{\ell,1}=\Sp_{2g_1}(\F_{\ell})$, alors $K_{\ell}^{\cap}=K(\mu_{\ell})$ donc on a $H_{\ell,2}=\Sp_{2g_2}(\F_{\ell})$. Mais par hypoth\`ese, ceci est impossible si $\ell$ est suffisamment grand. \findemo

\medskip

\noindent Nous introduisons maintenant un lemme qui nous servira dans la preuve du th\'eor\`eme \ref{isogenies1} ci-dessous et que nous r\'e-utilise\-rons dans le paragraphe 4. 

\begin{lemme}\label{reduction} Supposons qu'il existe $c>0$ et un ensemble $\Lambda$ infini tel que pour tout $\ell\in \Lambda$, on a 
\[ \text{Card}(H_{\ell,2})\leq c.\]
\noindent Alors, $g_1=g_2$ et pour tout $\ell\in\Lambda$ assez grand, on a 
\[ \ \textit{ soit }\ K_{\ell,1}=K_{\ell,2},\ \textit{ soit pour tout}\ i\in\{1,2\},\textit{ on a } \left[K_{\ell}^{\cup}:K_{\ell,i}\right]=2=\left[K_{\ell,i}:K_{\ell}^{\cap}\right].\]
\end{lemme}
\demo La preuve est une adaptation formelle de la remarque 2.4 de \cite{FJ}. Soit $\ell \in\Lambda$, $\ell\gg 1$. Consid\'erons  \`a nouveau le diagramme d'extensions pr\'ec\'edent. Par le lemme \ref{reduction0} pr\'ec\'edent, on sait que les groupes $H_{\ell,i}$ sont soit $\{I_{2g_i}\}$, soit $\{\pm I_{2g_i}\}$ (pour $i\in\{1,2\}$). On a 
\[\left|\Sp_{2g_1}(\F_{\ell})/H_{\ell,1}\right|=[K_{\ell}^{\cap}:K_{\mu_{\ell}}]=\left|\Sp_{2g_2}(\F_{\ell})/H_{\ell,2}\right|\]
\noindent et le lemme \ref{cardinal} implique donc que $g_1=g_2$ et que $H_{\ell,1}=H_{\ell,2}$. Le r\'esultat suit imm\'ediatement. \findemo

\medskip

\noindent Pour $i\in\{1,2\}$, nous noterons dans la suite $\rho_{\ell,i}$ la repr\'esentation $\ell$-adique modulo $\ell$ associ\'ee \`a la vari\'et\'e ab\'elienne $A_i$ et nous noterons $G_{\ell,i}$ son image dans $\textnormal{Aut}(A_i[\ell])$. 

\begin{theo}\label{isogenies1} Soient $A_1,A_2/K$ deux vari\'et\'es ab\'eliennes pleinement de type $\GSp$, sur un corps de nombres $K$. Soit $c>0$ telle qu'il existe un ensemble infini $\Lambda$ de nombres premiers, v\'erifiant 
\[\forall \ell\in \Lambda,\ \ \ \text{Card}(H_{\ell,2})\leq c.\]
\noindent Alors $A_1$ est $\bar K$-isog\`ene \`a $A_2$.
\end{theo}

\begin{rem} Dans le cas de dimension $1$, ceci est un r\'esultat de Frey-Jarden \cite{FJ} bas\'e sur les travaux de Serre \cite{serreim72}. En fait le r\'esultat de \cite{FJ} est plus g\'en\'eral car il vaut pour des courbes elliptiques quelconques sur un corps $K$ de type fini sur son sous-corps premier.
\end{rem}

\noindent Notons que l'on ne suppose pas a priori que $A_1$ et $A_2$ sont de m\^eme dimension : c'est une cons\'equence automatique.

\medskip

\noindent \demo Rappelons l'\'enonc\'e prouv\'e dans \cite{hrJIMJ} th\'eor\`eme 1.6 et remarque 1.7 que nous voulons utiliser : si $A_1$ et $A_2$ sont pleinement de type $\GSp$ de dimensions respectives $g_1$ et $g_2$ et non isog\`enes sur $\bar K$, alors pour tout $\ell$ suffisamment grand (d\'ependant de $A_1$ et $A_2$), on a
\[\rho_{\ell^{\infty},A_1\times A_2}\left(\Gal(K(A_1\times A_2[\ell^{\infty}])/K(\mu_{\ell^{\infty}}))\right)=\Sp_{2g_1}(\Z_{\ell})\times\Sp_{2g_2}(\Z_{\ell}).\]
\noindent Supposons donc par l'absurde que $A_1$ et $A_2$ ne sont pas $\bar K$-isog\`enes. Dans ce cas, pour tout $\ell$ assez grand, on a
\[\text{Card}\left(\rho_{\ell,A_1\times A_2}(G_K)\right)\gg\ll \ell^{\dim(\Sp_{2g_1})+\dim(\Sp_{2g_2})+1}=\ell^{2g_1^2+g_1+2g_2^2+g_2+1},\]
\noindent et de m\^eme, $A_1$ \'etant pleinement de type $\GSp$,
\[\text{Card}\left(\rho_{\ell,1}(G_K)\right)\gg\ll \ell^{2g^2_1+g_1+1}.\]
\noindent En particulier, on en d\'eduit que
\[\forall \ell\in \Lambda,\ \ \ c\geq \text{Card}(H_{\ell,2})=\left[K_{\ell}^{\cup}:K_{\ell,1}\right]\gg\ll \ell^{2g^2_2+g_2}.\]
\noindent Le dernier terme tendant vers l'infini avec $\ell$ ceci est contradictoire, donc $A_1$ et $A_2$ doivent \^etre $\bar K$-isog\`enes.\findemo

\medskip

\noindent Passons maintenant \`a la preuve alternative (B). \'Etant donn\'es un entier $n$ et deux repr\'e\-sen\-ta\-tions $\ell$-adiques modulo $\ell$, not\'ees $\rho_1,\rho_2 : G_K\rightarrow \GL_{n}(\F_{\ell})$, nous \'ecrirons $\rho_1\sim \rho_2$ s'il existe $u\in \GL_n(\F_{\ell})$ tel que $u^{-1}\rho_1u=\rho_2$. 

\begin{prop}\label{twist} Soit $\ell\geq 5$ un nombre premier tel que $G_{\ell,i}=\GSp_{2g_i}(\F_{\ell})$ pour $i\in\{1,2\}$. On suppose de plus que $g_1=g_2$ et que $[K_{\ell}^{\cup}:K_{\ell,1}]\leq 2$. Alors il existe un caract\`ere quadratique 
\[\epsilon_{\ell} : G_K\rightarrow \{\pm 1\}\ \text{ tel que }\ \rho_{\ell,1}\sim \epsilon_{\ell}\otimes\rho_{\ell,2}.\]
\noindent Si de plus $\ell\geq 4g+1$, alors ce caract\`ere $\epsilon_{\ell}$ est non ramifi\'e en toute place ultram\'etrique $v$ en laquelle $A_1$ et $A_2$ ont bonne r\'eduction et telle que $v$ est non ramifi\'ee sur $\Q$. En particulier, lorsque $\ell$ varie, le noyau $\ker(\epsilon_{\ell})$ varie dans un ensemble fini.
\end{prop}
\demo Pour la premi\`ere partie de l'\'enonc\'e, il s'agit d'une adaptation en dimension sup\'erieure du lemma 2.5 de \cite{FJ}, lui m\^eme bas\'e sur la preuve du lemme 8 de \cite{serreim72}. Pour la seconde partie de l'\'enonc\'e, concernant le caract\`ere quadratique $\epsilon_{\ell}$, il s'agit \'egalement de reprendre l'argument de la preuve du lemme 8 de \cite{serreim72} en utilisant le corollaire 3.4.4. de \cite{raynaud} en lieu et place des corollaires 11 et 12 de \cite{serreim72}. 

\noindent Nous utiliserons ici librement les notations introduites pr\'ec\'edemment et nous noterons $g=g_1=g_2$. Nous introduisons par ailleurs les projections  canoniques $\pi_{\ell,i} : \GSp_{2g_i}(\F_{\ell})\rightarrow \PGSp_{2g_i}(\F_{\ell})$ de noyau $\F_{\ell}^{\times}$ et notons $\overline{\rho}_{\ell,i}$ la compos\'ee de $\rho_{\ell,i}$ par $\pi_{\ell,i}$. Le corps fix\'e par $\ker(\overline{\rho}_{\ell,i})$ est not\'e $L_{\ell,i}$. Ce corps $L_{\ell,i}$ n'est autre que le sous-corps de $K_{\ell,i}$ fix\'e par le centre de $\Gal(K_{\ell,i}/K)$. Nous allons montrer que 
\begin{equation}\label{l1l2}
L_{\ell,1}=L_{\ell,1}.
\end{equation}
\noindent Admettons (\ref{l1l2}) un instant et voyons comment conclure \`a l'existence du caract\`ere quadratique $\epsilon_{\ell}$~: par la th\'eorie de Galois on obtient $\ker(\overline{\rho}_{\ell,1})=\ker(\overline{\rho}_{\ell,2})$. Ainsi, \`a composition \`a gauche par un automorphisme de $\PGSp_{2g}(\F_{\ell})$ pr\`es, $\overline{\rho}_{\ell,1}$ et $\overline{\rho}_{\ell,2}$ co\"{i}ncident. Or $\ell\geq 5$, donc d'apr\`es le lemme \ref{gsimple} les automorphismes de $\PGSp_{2g}(\F_{\ell})$ sont int\'erieurs et l'action par conjugaison se fait via un \'el\'ement de $\Sp_{2g}(\F_{\ell})$. Autrement dit, 
\[\exists y\in\Sp_{2g}(\F_{\ell})\ \forall x\in G_K,\ \  \pi_{\ell}\left(\rho_{\ell,1}(x)\right)=\pi_{\ell}\left(y^{-1}\rho_{\ell,2}(x)y\right).\]
\noindent Ainsi il existe un morphisme $\epsilon_{\ell} : G_K\rightarrow \F_{\ell}^{\times}$  tel que 
\[\forall x\in G_K,\ \ \rho_{\ell,1}(x)=\epsilon_{\ell}(x)\left(y^{-1}\rho_{\ell,2}(x)y\right).\]
\noindent En composant par le morphisme multiplicateur et en appliquant le lemme \ref{cyclo} on en d\'eduit que
\[\forall x\in G_K, \ \ \epsilon_{\ell}(x)^2=1.\]
\noindent Reste \`a prouver l'assertion (\ref{l1l2}) pour conclure la premi\`ere partie de la proposition. C'est ce que nous allons nous attacher \`a faire ci-dessous.

\medskip

\noindent Commen\c{c}ons par remarquer que si $K_{\ell,1}=K_{\ell,2}$, il n'y a rien \`a montrer. En utilisant les hypoth\`eses et le lemme \ref{reduction}, on sait donc que
\begin{equation} \label{sym} [K_{\ell,1}:K_{\ell}^{\cap}]=[K_{\ell,2}:K_{\ell}^{\cap}]=2.
\end{equation}

\noindent Notamment, la repr\'esentation $\rho_{\ell,1}$ envoie $K_{\ell}^{\cap}$ sur l'unique sous-groupe $\{\pm I_{2g}\}$ distingu\'e d'ordre $2$ de $\Sp_{2g}(\F_{\ell})$.

\medskip

\noindent \textbf{Fait 1 : } Le corps $K(\mu_{\ell})\cap L_{\ell,1}$ est l'unique sous-extension quadratique de $K(\mu_{\ell})/K.$

\medskip

\noindent Pour $\ell$ assez grand, l'extension $K(\mu_{\ell})/K$ est cyclique de groupe de Galois $\F_{\ell}^{\times}$. Elle poss\`ede donc une et une seule sous-extension quadratique correspondant \`a l'unique sous-groupe d'indice $2$ de $\F_{\ell}^{\times}$. On a le diagramme 
\[\xymatrix{
 & K_{\ell,1}\ar@{-}[dl]_{\Sp_{2g}(\F_{\ell})} \ar@{-}[dr]^{\F_{\ell}^{\times}=Z(\GSp_{2g}(\F_{\ell}))} & \\
K(\mu_{\ell})& & L_{\ell,1}}
 \]
\noindent d'extensions de corps. Or le groupe engendr\'e par $\Sp_{2g}(\F_{\ell})$ et $\F_{\ell}^{\times}$ (correspondant via la th\'eorie de Galois \`a l'extension $K(\mu_{\ell})\cap L_{\ell,1}$) n'est autre que le sous-groupe 

\[\{xM\in\GSp_{2g}(\F_{\ell})\ |\ x\in\F_{\ell}^{\times},\ M\in\Sp_{2g}(\F_{\ell}) \}.\] 

\noindent Ce sous-groupe est d'indice $2$ dans $\GSp_{2g}(\F_{\ell})$. En effet, il est clair sur la d\'efinition du multiplicateur $\lambda$ que $\lambda(xI_{2g})=x^2$ pour tout scalaire non nul $x$. D\`es lors on conclut en consid\'erant la d\'ecomposition $M^2=\lambda(M)\left(\frac{1}{\lambda(M)}M^2\right)$ ainsi que le fait que les carr\'es de $\F_{\ell}^{\times}$ forment un sous-groupe d'indice deux dans $\F_{\ell}^{\times}$ pour $\ell$ impair.\hfill$\Box$

\medskip

\noindent \textbf{Fait 2 : }On a l'\'egalit\'e $K_{\ell}^{\cap}=L_{\ell,1}\cdot K(\mu_{\ell}).$

\medskip

\noindent Par construction, $L_{\ell,1}$ est le sous-corps de $K_{\ell,1}$ des invariants par $\F_{\ell}^{\times}=Z(\GSp_{2g}(\F_{\ell})).$ Or on a montr\'e juste avant le fait 1 que $K_{\ell}^{\cap}$ est le sous-corps de $K_{\ell,1}$ des invariants par $\{\pm I_{2g}\}=Z(\Sp_{2g}(\F_{\ell}))$. On en d\'eduit que $L_{\ell,1}$ est inclus dans $K_{\ell}^{\cap}$. L'extension $K_{\ell,1}/K_{\ell}^{\cap}$ est de degr\'e deux et on voit en consid\'erant le diagramme d'extensions suivant
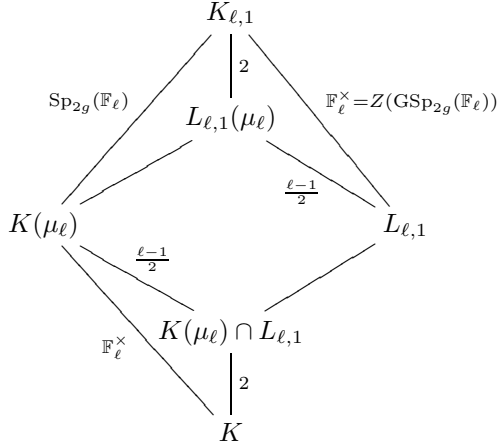
\[\xymatrix{
								 	& K_{\ell,1}\ar@{-}[ddl]_{\Sp_{2g}(\F_{\ell})} \ar@{-}[ddr]^{\F_{\ell}^{\times}=Z(\GSp_{2g}(\F_{\ell}))} \ar@{-}[d]^{2}&								 \\
 									& L_{\ell,1}(\mu_{\ell}) \ar@{-}[dl] \ar@{-}[dr]_{\frac{\ell-1}{2}}																												 &								 \\
K(\mu_{\ell})	\ar@{-}[ddr]_{\F_{\ell}^{\times}}\ar@{-}[dr]^{\frac{\ell-1}{2}}		&																							 		& \ar@{-}[dl] L_{\ell,1}	 \\
									& K(\mu_{\ell})\cap L_{\ell,1}	\ar@{-}[d]^{2}  																															&							 	 \\
									& K																																																							& }
 \]
\noindent que l'extension $K_{\ell,1}/L_{\ell,1}(\mu_{\ell})$ est \'egalement de degr\'e 2. Ceci conclut la preuve du fait 2.\hfill $\Box$

\medskip

\noindent L'\'equation (\ref{sym}) \'etant sym\'etrique en $A_1$ et $A_2$, on en d\'eduit sym\'etriquement que le corps $L_{\ell,2}$ v\'erifie \'egalement les faits $1$ et $2$. Le groupe de Galois $\Gal(K_{\ell}^{\cap}/K(\mu_{\ell}))=\PSp_{2g}(\F_{\ell})$ \'etant simple non ab\'elien, ceci d\'etermine le corps $L_{\ell,i}$ de mani\`ere unique et conclut la preuve de la premi\`ere partie de la proposition.

\medskip

\noindent Passons maintenant \`a la preuve de la seconde partie de l'\'enonc\'e. Soit $v$ une place ultram\'etrique du corps $K$ telle que $A_1$ et $A_2$ ont bonne r\'eduction en $v$ et que $v$ est non ramifi\'ee sur $\Q$ (ie telle que si $p$ est un premier tel que $v|p$ alors l'indice de ramification $e(v/p)$ est \'egal \`a $1$). Supposons que la caract\'eristique de $v$ est $\ell$ (en effet $\epsilon_{\ell}$ est non ramifi\'e en $v$ sinon par le theorem 1 de \cite{serretate}) et notons $\bar{\F}_{\ell}$ une cl\^oture alg\'ebrique de $\F_{\ell}$. Suivant \cite{serreim72} paragraphe 1.13. et \cite{raynaud} corollaire 3.4.4. nous notons par ailleurs $\chi_1,\ldots,\chi_{2g}$ (respectivement $\chi_1',\ldots,\chi_{2g}'$) les caract\`eres du groupe d'inertie mod\'er\'ee en $v$ \`a valeurs dans $\bar{\F}_{\ell}^{\times}$, intervenant dans le module galoisien $A_1[\ell]\otimes \bar{\F}_{\ell}$ (resp. $A_2[
 \ell]\otimes \bar{\F}_{\ell}$),  En notant $\epsilon_v$ la restriction de $\epsilon_{\ell}$ au groupe d'inertie en $v$, on a pour tout $i$, (quitte \`a renum\'eroter les $\chi_i'$)
\[\chi_i=\epsilon_v\chi_i'.\]
\noindent Comme l'indice de ramification de $v$ est $1$, le corollaire 3.4.4. de \cite{raynaud} nous dit que les $\chi_i$ sont de la forme
\[\chi_i=\theta_{k_1}^{e(k_1)}\ldots\theta_{k_n}^{e(k_n)}\]
\noindent o\`u pour tout $r$, $e(r)\in\{0,1\}$ et o\`u les $\theta_{k_i}$ sont les $n$ caract\`eres fondamentaux de niveau $n$ (cf. \cite{serreim72} p.267 et \cite{raynaud} d\'efinition 1.1.1), l'entier $n$ pouvant varier dans $1,\ldots,2g$. Les invariants des $\chi_i$ et $\chi_i'$ dans $\Q/\Z$ (cf. \cite{serreim72} paragraphe 1.7) varient donc dans l'ensemble :
\[X=\left\{e(k_1)\frac{\ell^{k_1}}{\ell^n-1}+\cdots+e(k_n)\frac{\ell^{k_n}}{\ell^n-1}\ |\ k_i\in\{0,\ldots,n-1\},\ n\in\{1,\ldots,2g\}\right\}.\]
\noindent Enfin, l'invariant de $\epsilon_v$, $\text{Inv}(\epsilon_v)$ vaut $0$ ou $\frac{1}{2}$  car $\epsilon_v^2=1$. Par ailleurs, $\text{Inv}(\epsilon_v)$ est de la forme $x-x'$ avec $x=\text{Inv}(\chi_i)$ et $x'=\text{Inv}(\chi_i')$. En particulier on a $x,x'\in X$. Or si 
\[x=e(k_1)\frac{\ell^{k_1}}{\ell^n-1}+\cdots+e(k_n)\frac{\ell^{k_n}}{\ell^n-1}, \text{ on a alors}\]
\noindent 
\[0\leq x\leq \frac{n\ell^{n-1}}{\ell^n-1}<\frac{2g}{\ell-1}.\]
\noindent En particulier, $|x-x'|<\frac{2g}{\ell-1}$ et comme $\ell\geq 4g+1$, on voit que $|x-x'|<\frac{1}{2}$, donc n\'ecessairement l'invariant de $\epsilon_v$ vaut $0$ ce qui signifie que $\epsilon_v$ est non ramifi\'e en $v$.

\medskip

\noindent Pour la derni\`ere assertion de la proposition, notons $S$ l'ensemble des places de bonne r\'eduction de $A_1$ et $A_2$. Nous venons de prouver que les caract\`eres $\epsilon_{\ell}$ se factorisent \`a travers le groupe de Galois $G_S$ de la sous-extension maximale de $\bar K$ non ramifi\'ee en dehors de $S$. Mais par le th\'eor\`eme classique de Hermite, il n'existe qu'un nombre fini d'extensions de degr\'e au plus $2$ de $K$, non ramifi\'ees en dehors de $S$. Autrement dit les noyaux de $\epsilon_{\ell}$ varient dans un ensemble fini quand $\ell$ est variable.\findemo

\medskip

\noindent Nous utiliserons ci-dessous la version suivante des r\'esultats de Faltings \cite{falt} sur la conjecture de Tate pour les vari\'et\'es ab\'eliennes. 

\begin{prop}\label{lemfalt} Soit $K$ un corps de nombres et soient $A/K$ et $B/K$ deux vari\'et\'es ab\'eliennes.
\begin{itemize}
\item Si $\rho_{\ell^{\infty},A}$ est isomorphe \`a $\rho_{\ell^{\infty},B}$ alors $A$ est $K$-isog\`ene \`a $B$.
\item Il existe $C_0=C_0(A,K)$ telle que si pour un nombre premier $\ell\geq C_0$ les repr\'esentations modulo $\ell$, $\rho_{\ell,A}$ et $\rho_{\ell,B}$ sont isomorphes, alors $A$ est $K$-isog\`ene \`a $B$.
\end{itemize}
\end{prop}
\demo Dans l'article de Faltings \cite{falt}, le premier \'enonc\'e est d\'emontr\'e dans le Korollar 2, page 361~; le deuxi\`eme \'enonc\'e, pour les repr\'esentations modulo $\ell$,  peut se d\'eduire des d\'emonstrations comme cela est montr\'e par Zarhin \cite{zar1}, Corollary 5.4.5. 
\findemo

\medskip

\noindent \textbf{Preuve (B) du th\'eor\`eme 3.1} Soit $\Lambda$ un ensemble infini de premiers comme dans l'\'enonc\'e et soit $\ell\in \Lambda$, $\ell\geq 4g+1$. D'apr\`es le lemme de r\'eduction \ref{reduction}, on voit que $A_1$ et $A_2$ ont m\^eme dimension et que pour $i\in\{1,2\}$, on a $[K_{\ell}^{\cup}: K_{\ell,i}]\in\{1,2\}$ (quitte \`a enlever un nombre fini de premiers de l'ensemble $\Lambda$). La proposition \ref{twist} entra\^ine alors qu'il existe un caract\`ere quadratique 
\[\epsilon_{\ell} : G_K\rightarrow \{\pm 1\}\ \text{ tel que }\ \rho_{\ell,1}\sim \epsilon_{\ell}\otimes\rho_{\ell,2},\]
\noindent et que de plus quand $\ell$ est variable, les $\epsilon_{\ell}$ varient dans un ensemble fini. Quitte \`a remplacer $\Lambda$ par une sous-partie infinie, on peut donc supposer que le caract\`ere $\epsilon_{\ell}$ est ind\'ependant de $\ell$~: notons le $\epsilon$. Sur l'extension $K'$ de $K$ correspondant au noyau de $\epsilon$, on obtient donc que les repr\'esentations $\rho_{\ell,1}$ et $\rho_{\ell,2}$ sont isomorphes. Par la seconde partie de la proposition \ref{lemfalt} ceci implique que $A_1$ et $A_2$ sont $K'$-isog\`enes, ce qui conclut la d\'emonstration.\hfill $\Box$

\section{Preuve du th\'eor\`eme \ref{radical}}
\noindent Nous conservons dans la suite les notations de la partie pr\'ec\'edente (notamment $A_1$ et $A_2$ sont pleinement de type GSp), nous notons $G_{\ell,A_1\times A_2}\subset G_{\ell,1}\times G_{\ell,2}$ le groupe de Galois correspondant \`a $A_1\times A_2$ et $\pi_i$ la projection de $G_{\ell,A_1\times A_2}$ sur $G_{\ell,i}$ pour $i\in\{1,2\}$. Rappelons que comme indiqu\'e dans le paragraphe \ref{structure} concernant la structure de la preuve (de la partie suffisante) du th\'eor\`eme principal \ref{radical}, nous d\'ecomposons cette preuve en trois \'etapes. Dans la suite, on suppose donn\'e un sous-ensemble $S$ de places finies $v$ de $K$, de corps r\'esiduel $\F_v$, de bonne r\'eduction pour $A_1$ et $A_2$, de densit\'e analytique $1$ et nous supposons \'egalement donn\'e un sous-ensemble $\Lambda$ infini de l'ensemble des nombres premiers tel que
\begin{equation}\label{hypo}\forall v\in S,\ \forall \ell\in\Lambda\ \ \left( \ell\mid \text{Card}(A_1(\F_v))\iff \ell\mid \text{Card}(A_2(\F_v))\right).
\end{equation}
\noindent Nous utiliserons le lemma 4.6 de \cite{hallp} que nous rappelons ci-dessous.

\begin{lemme} \textnormal{\textbf{(\cite{hallp})}} \label{lemmehp} Soit $\ell\in\Lambda$. Sous l'hypoth\`ese (\ref{hypo}) ci-dessus, on a 
 \[\forall (x,y)\in G_{\ell,A_1\times A_2},\ \ \det(x-I_{2g_1})=0\iff\det(y-I_{2g_2})=0.\]
\end{lemme}

\subsection{Preuve de l'\'etape (\ref{1})}

\begin{lemme} \label{fait1} Soit $\ell\geq 5$ premier ne divisant pas le degr\'e des polarisations fix\'ees sur $A_1$ et $A_2$, tel que $K_{\ell,1}\not\subset K_{\ell,2}$. Le groupe $\pi_1(\ker(\pi_2))$ est distingu\'e dans $G_{\ell,1}$ et contient l'\'el\'ement non trivial $-I_{2g_1}$.
\end{lemme}
\demo La projection $\pi_1$ est surjective et $\ker(\pi_2)$ est distingu\'e dans $G_{\ell,A_1\times A_2}$, ce qui entra\^ine la premi\`ere partie de l'asser\-tion. Si $u\in\pi_1(\ker(\pi_2))$, il existe $x\in G_K$ tel que $u=\rho_{\ell,1}(x)$ et $\rho_{\ell,2}(x)=I_{2g_2}$. De plus on sait par le lemme \ref{cyclo} que 
$\lambda\circ\rho_{\ell,1}=\lambda\circ\rho_{\ell,2}$, donc en particulier, $\lambda(u)=\lambda(I_{2g_2})=1$. Ceci prouve que $\pi_1(\ker(\pi_2))$ est inclus dans le groupe $\Sp_{2g_1}(\F_{\ell})$.  Montrons maintenant que $\pi_1(\ker(\pi_2))$ est non trivial. En effet, dire que $\pi_1(\ker(\pi_2))=\{I_{2g_1}\}$ \'equivaut \`a dire que $\ker(\pi_2)\subset \ker(\pi_1)$ ce qui \'equivaut visiblement \`a dire que $\ker(\pi_2)$ est trivial. Dans ce cas le diagramme suivant est commutatif
\[\xymatrix{
G_K\ar[dr]_{\rho_{\ell,2}} \ar[rr]^{\rho_{\ell,1}} & & G_{\ell,1} \\
& G_{\ell,2}\ar[ur]_{\pi_1\circ\pi_2^{-1}} & }
 \]
\noindent et on en d\'eduit que $\ker(\rho_{\ell,2})\subset \ker(\rho_{\ell,1})$. Ceci entra\^ine que l'extension $K_{\ell,1}$ est incluse dans $K_{\ell,2}$ ce qui contredit l'hypoth\`ese. Ainsi comme $\ell\geq 5$, le lemme \ref{gsimple} implique que le groupe $\pi_1(\ker(\pi_2))$ est $\{\pm I_{2g_1}\}$ ou $\Sp_{2g_1}(\F_{\ell})$ et notamment $-I_{2g_1}\in\pi_1(\ker(\pi_2)).$\findemo

\medskip

\noindent Nous pouvons maintenant passer \`a la preuve de l'\'etape (1) proprement dite~: supposons par l'absurde que pour tout $\ell\in\Lambda$ assez grand on a $K_{\ell,1}\not= K_{\ell,2}$. On peut notamment supposer que $\ell\geq 5$ et quitte \`a intervertir $A_1$ et $A_2$, le lemme \ref{fait1} nous assure que l'\'el\'ement $(-I_{2g_1},I_{2g_2})$ est dans $G_{\ell,A_1\times A_2}$. Mais $\det(-I_{2g_1}-I_{2g_1})\not=0$ alors que $\det(I_{2g_2}-I_{2g_2})=0$. Ceci contredit le lemme \ref{lemmehp} pr\'ec\'edent et conclut la preuve de la premi\`ere partie de l'\'etape (\ref{1}). Pour une infinit\'e de $\ell$ il est donc vrai que $K_{\ell,1}=K_{\ell,2}$. Notre th\'eor\`eme \ref{isogenies1} sur les isog\'enies horizontales permet de conclure que $A_1$ et $A_2$ sont $\bar{K}$-isog\`enes (en effet pour les $\ell$ concern\'es on a $[K_{\ell}^{\cup}:K_{\ell,1}]=1$). Nous utilisons ici de mani\`ere sous-jacente les r\'esultats de Faltings \cite{falt} (proposition \ref{lemfalt} pr\'ec\'edente).\hfill $\Box$

\subsection{Preuve de l'\'etape (\ref{2})}

\noindent Nous savons maintenant que les vari\'et\'es ab\'eliennes $A_1$ et $A_2$ sont isog\`enes sur $\bar{K}$. En particulier elles ont m\^eme dimension, not\'ee $g$. De plus par l'\'etape (\ref{1}) on sait qu'il existe un sous-ensemble infini $\Lambda_1$ de $\Lambda$ tel que pour $\ell\in\Lambda_1$ on a $[K_{\ell}^{\cup}:K_{\ell,1}]=1$. Pour un tel premier $\ell\geq 5$ la proposition \ref{twist} affirme l'existence d'un caract\`ere quadratique 
\[\epsilon_{\ell} : G_K\rightarrow \{\pm 1\}\ \text{ tel que }\ \rho_{\ell,1}\sim \epsilon_{\ell}\otimes\rho_{\ell,2}.\]
\noindent De plus en choisissant $\ell\geq 4g+1$, cette m\^eme proposition assure que les caract\`eres $\epsilon_{\ell}$ varient dans un ensemble fini (quand $\ell$ est variable dans $\Lambda_1$). Quitte \`a remplacer $\Lambda_1$ par une sous-partie infinie, on peut donc supposer que $\epsilon_{\ell}$ est ind\'ependant de $\ell$. C'est pr\'ecis\'ement ce que l'on voulait montrer. \hfill $\Box$

\subsection{Preuve de l'\'etape (\ref{3})}
\noindent Montrons que $\rho_{\ell,1}$ est isomorphe \`a $\rho_{\ell,2}$ pour une infinit\'e de $\ell$. La conclusion r\'esultera alors de la proposition \ref{lemfalt}. Pour ce faire nous utiliserons le r\'esultat suivant sur la d\'ecomposition d'une place dans une extension de la forme $K(A[\ell])$.

\begin{lemme} \label{td} Soient $p\not=q$ deux nombres premiers, soit $K$ un corps de nombres et soit $A/K$ une vari\'et\'e ab\'elienne. Soit $v$ une place finie de $K$ au dessus de $p$, de corps r\'esiduel $\F_v$, telle que $A$ a bonne r\'eduction en $v$. Le groupe $A(\F_{v})$ contient $A[q]\simeq (\Z/q\Z)^{2\dim A}$ si et seulement si $v$ est totalement d\'ecompos\'ee dans l'extension $K(A[q])/K$.
\end{lemme}
\demo Notons $\pi_{v}$ l'endomorphisme de Frobenius sur la r\'eduction $A_{v}/\F_{v}$ de $A/K$ en une place de bonne r\'eduction $v$. Notons que $p$ \'etant premier \`a $q$, la r\'eduction modulo $v$ est injective sur les points de $q$-torsion et on identifie $A_{v}[q]$ et $A[q]$ via cette injection. Par construction on a 
\[\ker(\pi_{v}-\text{Id})=A_{v}(\F_{v}).\]
\noindent Donc $A_{v}(\F_{v})$ contient le groupe $A_v[q]$ si et seulement si $\pi_{v | A_{v}[q]}=\text{Id}_{A_{v}[q]}$. Cette derni\`ere assertion est \'equivalente \`a dire que le groupe de d\'ecomposition $D_q(v/p)$ est trivial dans $\Gal(K(A[q])/K)$ et ceci \'equivaut \`a dire que  $v$ est totalement d\'ecompos\'ee dans l'extension $K(A[q])/K$. \findemo

\begin{lemme} \label{galois} Soit $A/K$ une vari\'et\'e ab\'elienne pleinement de type $\GSp$ et $L/K$ une extension finie. Si $\ell$ est suffisamment grand alors $L\cap K(A[\ell])=K$.
\end{lemme}
\demo Il s'agit de la proposition 2.6 de \cite{hallp} dans laquelle il faut appliquer le corollaire du th\'eor\`eme 3 de \cite{college8586} en lieu et place du th\'eor\`eme 3 de \cite{serreim72}.\findemo

\medskip

\noindent Pour une vari\'et\'e ab\'elienne $A$, notons $S_K(A)$ l'ensemble des places finies de $K$ de bonne r\'eduction pour $A$ et notons 
\[ S_{\ell}^{td}:=\left\{v\in S_K(A_1)\cap S_K(A_2)\ |\ v\text{ se d\'ecompose totalement dans }K_{\ell,1}\right\}.\]
\noindent D'apr\`es le lemme \ref{td} pr\'ec\'edent, si $v\in S_{\ell}^{td}$ alors $A_1(\F_v)[\ell]=A_1[\ell]$. En particulier on voit que $\ell\mid\text{Card}(A_1(\F_v))$. Pour $i\in\{1,2\}$ et en notant $v_{\ell}$ la valuation $\ell$-adique sur les entiers, introduisons
\[f_{\ell,i} : S_K(A_i)\rightarrow \{0,1\},\ \text{ d\'efinie par }\ v\mapsto \min\left\{1,v_{\ell}(\text{Card}(A_i(\F_v)))\right\}.\]
\noindent Par construction $f_{\ell,i}(v)=1$ si et seulement si $\ell\mid \text{Card}(A_i(\F_v))$. Nous venons donc de prouver que si $v\in S_{\ell}^{td}$ alors $f_{\ell,1}(v)=1$.

\medskip

\noindent Montrons maintenant que nos hypoth\`eses impliquent que $f_{\ell,2}(v)=0$ pour $v$ variant dans un sous-ensemble de densit\'e strictement positive, ce qui contredira l'hypoth\`ese (\ref{hypo}). Pr\'ecis\'ement, on sait qu'il existe un caract\`ere $\epsilon : G_K\rightarrow \{\pm 1\}$ tel que $\rho_{\ell,1}\sim\epsilon\otimes\rho_{\ell,2}$. Supposons que $\epsilon$ est non  trivial, notons $L/K$ l'extension quadratique correspondant au caract\`ere $\epsilon$ et notons $\zeta_v:=\epsilon(\Frob_v)\in\{\pm 1\}$. Notons $S_{\ell}'$ le sous-ensemble de $S_{\ell}^{td}$ constitu\'e des $v$ telles que $\zeta_v=-1$. L'ensemble $S_{\ell}'$ a une densit\'e non nulle d'apr\`es le lemme \ref{galois}. Reste \`a v\'erifier que $f_{\ell,2}$ est nulle sur cet ensemble.

\medskip

\noindent La place $v$ est totalement d\'ecompos\'ee dans $K_{\ell,1}$ donc $\rho_{\ell,1}(\Frob_v)=\Id_{A_1[\ell]}$. De plus le Frobenius est un g\'en\'erateur topologique du groupe de Galois $G_{\F_v}$ et via l'accouplement de Weil on sait que le groupe $\mu_{\ell}$ est inclus dans $A_1[\ell]$. Ceci nous permet d'en d\'eduire que si $\sigma\in G_{\F_v}$ alors $\sigma$ agit trivialement sur $\mu_{\ell}$. Notons par ailleurs 
\[\pi_{A_i}(T)=\det\left(\text{Id}-T\rho_{\ell^{\infty},i}(\Frob_v)|V_{\ell}(A_i)\right)\ \text{ pour }\ i\in\{1,2\}.\]
\noindent On voit que $v$ totalement d\'ecompos\'ee implique que 
\[\pi_{A_1}(T)\equiv \det(\text{Id}-T\text{Id})\equiv (1-T)^{2g}\mod\ell.\]
\noindent Rappelons que $\rho_{\ell,1}(\Frob_v)$ est, \`a conjugaison pr\`es, \'egal \`a $\epsilon(\Frob_v)\rho_{\ell,2}(\Frob_v)$. Notamment, si de plus $v\in S_{\ell}'$ on a $\epsilon(\Frob_v)=-1$ et donc 
\[ \pi_{A_2}(T)\equiv (1+T)^{2g}\mod\ell.\]
\noindent En particulier, 
\[\text{Card}(A_2(\F_v))\equiv\pi_{A_2}(1)\equiv 2^{2g}\not\equiv 0\mod\ell.\]
\noindent Ceci prouve que $f_{\ell,2}(v)=0$ ce qui conclut. \hfill$\Box$

\end{document}